\documentclass[11pt]{article}
\usepackage[letterpaper,top=0.89in,bottom=0.91in,left=0.96in,right=0.96in]{geometry}
\usepackage[T1]{fontenc}
\usepackage{lmodern,microtype,amsmath,amssymb,amsthm,mathtools}
\usepackage[hidelinks]{hyperref}
\usepackage{bookmark}
\allowdisplaybreaks[1]
\newcommand{\C}{\mathbb C}
\newcommand{\ket}[1]{\lvert #1\rangle}
\newcommand{\Span}{\operatorname{span}}
\newcommand{\wt}{\operatorname{wt}}
\newcommand{\gr}{\operatorname{gr}}
\newcommand{\Sym}{\operatorname{Sym}}
\newcommand{\Supp}{\operatorname{Supp}}
\newcommand{\Ann}{\operatorname{Ann}}
\newtheorem{theorem}{Theorem}[section]
\newtheorem{lemma}[theorem]{Lemma}
\newtheorem{proposition}[theorem]{Proposition}
\newtheorem{corollary}[theorem]{Corollary}
\newtheorem*{restatedquasilisse}{Corollary~\ref{cor:quasilisse}}
\newtheorem*{restatedrational}{Theorem~\ref{thm:rational}}
\newtheorem*{restatedtwisted}{Theorem~\ref{thm:twisted}}
\newtheorem*{restatedav}{Theorem~\ref{thm:av}}
\newtheorem*{restatedsfg}{Theorem~\ref{thm:sfg}}
\theoremstyle{remark}
\newtheorem{remark}[theorem]{Remark}
\hypersetup{pdftitle={Associated Varieties of Finite-Group Orbifolds},
 pdfauthor={Andrew Riesen}}
\title{Associated Varieties of Finite-Group Orbifolds}
\author{Andrew Riesen}
\date{\today}
\begin{document}
\maketitle

\begin{abstract}
Let $V$ be a simple, strongly finitely generated vertex algebra with a connected nonnegative integer grading, and let $G$ be a finite group of grading-preserving automorphisms. In this paper, we prove that $X_{V^G}\cong X_V/G$ as reduced Poisson varieties. More precisely, the canonical surjection
$
R_{V^G}\longrightarrow (R_V)^G $
has nilpotent kernel. The proof uses a uniform derivative bound in a relative $C_2$-quotient, which also implies that $V^G$ is strongly finitely generated. Consequently, quasi-lisseness and \(C_2\)-cofiniteness descend to finite-group orbifolds, and hence, by existing results, so does strong rationality.
\end{abstract}

\section{Introduction}\label{sec:intro}

A central theme in invariant theory is understanding which properties of an algebraic object are preserved under taking invariants. The setting of interest in this paper is that of a finite group $G$ acting on a vertex algebra $V$, an algebraic structure arising naturally in two-dimensional conformal field theory. The invariant subalgebra $V^G$, also sometimes called an \emph{orbifold} (although this differs from the meaning of the term in physics), gives rise to several fundamental questions. For example, if the representation theory of $V$ is semisimple and has only finitely many irreducible modules, does the same hold for $V^G$? Similarly, if $V$ is strongly finitely generated, must $V^G$ also be strongly finitely generated? The latter is a vertex-algebra analogue of Hilbert's classical theorem on finite generation of invariants. These and related questions have been the subject of substantial work over the years \cite{ADFLM,AL,CM,DLM,DM94,GR,HuangProblems,LianLinshaw,McRae,Miyamoto,Xu}. In this paper, we continue this work by describing the associated variety of $V^G$ in terms of the associated variety of $V$ and the action of $G$, as well as showing strong finite generation descends.  

We work throughout with ordinary vertex algebras over $\C$. By a connected nonnegatively graded vertex algebra, we mean a vertex algebra $V$ equipped with a grading
\[
V=\bigoplus_{n\geq 0}V_n
\]
such that $V_0=\C\ket{0}$ and $a_mV_k\subseteq V_{h+k-m-1}$ for $a\in V_h$, where $V_j=0$ for $j<0$. For a vertex algebra $W$, Zhu's $C_2$-algebra
\cite{Zhu} is defined by
\[
 C_2(W)=\Span_{\C}\{a_{-2}b\mid a,b\in W\},
 \qquad R_W=W/C_2(W).
\]
The algebra $R_W$ carries a natural commutative Poisson algebra
structure, with multiplication and Poisson bracket given by
\[
 \overline a\cdot\overline b=\overline{a_{-1}b},
 \qquad \{\overline a,\overline b\}=\overline{a_0b},
\]
where $\overline a$ denotes the image of $a$ in $R_W$
\cite{Li,ArakawaC2}. The associated variety of $W$ is the
reduced affine Poisson variety
\[
 X_W:=\operatorname{Spec}\bigl((R_W)_{\mathrm{red}}\bigr).
\]
For a finite group $G$ acting on $V$, write
$X_V/G=\operatorname{Spec}\bigl(((R_V)_{\mathrm{red}})^G\bigr)$
for the affine quotient of its associated variety. Our main result
is the following.

\begin{theorem}[Associated varieties of finite-group orbifolds]\label{thm:av}
Let $V$ be a simple, strongly finitely generated vertex algebra with
a connected nonnegative integer grading, and let $G$ be a finite
group of grading-preserving automorphisms. The canonical Poisson
homomorphism $p_G:R_{V^G}\twoheadrightarrow(R_V)^G$ has nilpotent
kernel. In particular it induces an isomorphism of reduced Poisson
varieties
\[
 X_V/G\ \xrightarrow{\ \sim\ }\ X_{V^G}.
\]
\end{theorem}

This proves Conjecture~13.2 of Lian and Linshaw
\cite{LianLinshaw}. The map $p_G$ is induced by the inclusion $V^G\subseteq V$,
and its surjectivity follows from the exactness of taking
finite-group invariants. Its kernel is given by $\ker p_G=\frac{V^G\cap C_2(V)}{C_2(V^G)}$, in other words the kernel records the additional $C_2$-relations obtained
by allowing fields from all of $V$, rather than only $V^G$.
The proof of Theorem~\ref{thm:av} uses a uniform bound on
derivatives in a relative $C_2$-quotient 
(Proposition~\ref{prop:cutoff}), namely $R_{V^G}(V)$ where $V$ is considered as a $V^G$-module. The same uniform bound also gives
descent of strong finite generation.

\paragraph{Descent of strong finite generation.}
A vertex algebra is strongly finitely generated if it is spanned by
normally ordered products of derivatives of finitely many fields.
As mentioned above, the descent of strong finite
generation is a vertex-algebra analogue of Hilbert's
theorem. Linshaw and collaborators have studied this
problem extensively, establishing strong finite
generation and related derivative-order bounds for
several important families of orbifolds by arbitrary complex
reductive groups, including free-field algebras,
affine vertex (super)algebras, generic-level
$\mathcal W$-(super)algebras, and their tensor products
\cite{LinshawHilbert,LinshawHeisenberg,LinshawAffine,
LinshawKWY,CreutzigLinshawSF,CreutzigLinshawTrialities}.
The following theorem establishes strong finite
generation uniformly for finite-group orbifolds
under our hypotheses.

\begin{theorem}[Finite strong generation]\label{thm:sfg}
Let $V$ be a simple, strongly finitely generated vertex algebra
with a connected nonnegative integer grading, and let $G$ be a
finite group of grading-preserving automorphisms. Then $V^G$
is strongly finitely generated.
\end{theorem}

Theorems~\ref{thm:av} and~\ref{thm:sfg} together have several
consequences for finite-group orbifolds. We begin with the
quasi-lisse property.

\paragraph{Quasi-lisse orbifolds and symplectic leaves.}

Building on the work of Zhu and Li
\cite{Zhu,Li,LiAbelianizing}, Arakawa introduced the
associated variety $X_V$ of a vertex algebra and gave
a geometric interpretation of $C_2$-cofiniteness
\cite{ArakawaC2}. Namely, for a strongly finitely
generated vertex algebra with a connected nonnegative
integer grading, $C_2$-cofiniteness is equivalent to
the associated variety being a point. Generalizing
this condition, Arakawa and Kawasetsu introduced the
notion of a quasi-lisse vertex algebra: a strongly
finitely generated vertex algebra whose associated
variety, viewed as an affine Poisson variety, has
finitely many symplectic leaves \cite{AK}.
They showed that a simple quasi-lisse VOA of CFT type
has only finitely many simple ordinary modules, and
that their normalized characters satisfy a common
modular linear differential equation \cite{AK}.
Following a suggestion of Dra\v{z}en Adamovi\'c,
Arakawa, Dai, Fasquel, Li, and Moreau conjectured
that the quasi-lisse property descends to orbifolds
by finite solvable groups \cite[Conjecture~B]{ADFLM}.

Since finite Poisson quotients preserve the finiteness
of symplectic leaves \cite{ES2010},
Theorems~\ref{thm:av} and~\ref{thm:sfg} together
establish descent of quasi-lisse in the finite group case.

\begin{corollary}[Quasi-lisse descent]
\label{cor:quasilisse}
Let $V$ be a simple quasi-lisse vertex algebra with
a connected nonnegative integer grading, and let
$G$ be a finite group of grading-preserving automorphisms.
Then $V^G$ is quasi-lisse.
\end{corollary}

\paragraph{Strong rationality and vertex tensor categories.}
A central question in the theory of vertex operator
algebras has been whether strong rationality descends
from a VOA $V$ to its fixed-point
subalgebra $V^G$.
Recall that a VOA is called \emph{strongly rational}
if it is simple, of CFT type, self-contragredient,
rational, and $C_2$-cofinite. By a theorem of Huang,
the category of ordinary modules of a strongly
rational VOA is a modular tensor category
\cite{HuangRigidity}. Thus, descent of strong
rationality would in particular imply that the
category of ordinary $V^G$-modules is modular.

For brevity, we refer to the general question of
descent of strong rationality as the
\emph{orbifold problem}.
Progress towards its resolution was made by Miyamoto,
who established $C_2$-cofiniteness for cyclic
groups, and consequently for solvable groups by iteration \cite{Miyamoto}.
Carnahan and Miyamoto subsequently proved regularity
for solvable-group orbifolds under suitable
hypotheses \cite{CM}. For arbitrary finite groups,
McRae showed that strong rationality of $V^G$
follows from that of $V$, provided that $V^G$
is $C_2$-cofinite \cite[Corollary~4.23]{McRae}.

Our description of the associated variety supplies
this remaining condition. Indeed, if $V$ is
$C_2$-cofinite, then $X_V$ is a point, and
Theorem~\ref{thm:av} implies that $X_{V^G}$ is
also a point. Since $V^G$ is strongly finitely generated by
Theorem~\ref{thm:sfg}, Lemma~\ref{lem:stronggen-lift}
shows that $R_{V^G}$ is finitely generated. As its
reduced spectrum is a point, $R_{V^G}$ is finite
dimensional, so $V^G$ is $C_2$-cofinite. Combined with
McRae's result, we obtain the following theorem,
resolving the orbifold problem.

\begin{theorem}[Descent of Strong Rationality]
\label{thm:rational}
Let $V$ be a strongly rational VOA and
$G$ a finite group of VOA automorphisms.
Then $V^G$ is strongly rational.
\end{theorem}

The representation theory of $V^G$ is also closely connected to another important feature of orbifold constructions: the appearance of twisted sectors. In the language of
vertex operator algebras, the corresponding objects are $g$-twisted
$V$-modules whose monodromy is determined by the action of $g\in G$. 
The mathematical theory of twisted modules has been developed
extensively, with foundational contributions including the work
of Frenkel, Lepowsky, and Meurman \cite{FLM} and Dong, Li, and
Mason \cite{DLMtwisted,DLM}. In particular, every $g$-twisted $V$-module restricts
naturally to a module for $V^G$, and it is natural to ask how
the tensor products of the different twisted sectors are related
to the representation theory of $V^G$.

Huang conjectured that, for a suitable VOA and finite group $G$,
the categories of $g$-twisted modules should together form a
$G$-crossed braided tensor category
\cite[Conjecture~6.2]{HuangProblems}. Du and Huang subsequently
developed the theory of twisted intertwining operators and their
tensor products \cite{DuHuang}. Recently, Xu constructed a
$G$-crossed braided fusion category of twisted modules for
strongly rational VOAs \cite{XuCrossed}, and established a
further crossed modularity result \cite{XuCrossedModularity}.

McRae's equivariantization theorem provides another approach
to the $G$-crossed tensor category structure, assuming that the
fixed-point algebra admits a suitable vertex tensor category
of modules \cite{McRae}. By our $C_2$-cofiniteness result,
together with Huang's logarithmic tensor product theory
\cite{HuangCofinite,HLZ}, these hypotheses are satisfied
under the assumptions below, without requiring rationality.

Write $\operatorname{Mod}(W)$ for the category of
grading-restricted generalized modules of a VOA $W$, and
$\operatorname{TwMod}_g(V)$ for the corresponding category of
$g$-twisted $V$-modules.

\begin{theorem}[Twisted sectors of a finite orbifold]\label{thm:twisted}
Let $V$ be a simple $C_2$-cofinite VOA of CFT type and
$G\leq\operatorname{Aut}(V)$ a finite group. Then
\[
 \operatorname{TwMod}_G(V)
 =\bigoplus_{g\in G}\operatorname{TwMod}_g(V)
\]
admits a natural $G$-crossed braided tensor category structure.
Its equivariantization is braided tensor equivalent to
$\operatorname{Mod}(V^G)$.
\end{theorem}

In the strongly rational setting, this gives an alternative proof
of the existence of the $G$-crossed braided tensor category
structure, using McRae's categorical approach. More generally,
Theorem~\ref{thm:twisted} applies without assuming rationality
or semisimplicity, and hence also applies in logarithmic
conformal field theory.

A further consequence concerns the reconstruction of twisted
group doubles by holomorphic orbifolds. Theorem~\ref{thm:rational}
removes the orbifold-rationality hypothesis in the reconstruction
theorem of Gannon and the author \cite[Theorem~3.2]{GR}, based on
Evans--Gannon \cite{EG}. We record the resulting unconditional
realization in Section~\ref{subsec:holomorphic}.

Finally, there is also a consequence for genus-one correlation functions.
The $C_2$-cofiniteness of $V^G$ allows us to apply Huang's modular
invariance theorem for logarithmic intertwining operators
\cite{HuangLog}. In particular, the genus-one correlation functions for $V^G$ constructed
from shifted pseudo-traces span spaces closed under modular
transformations (Corollary~\ref{cor:genusone}). Here pseudo-traces generalize
ordinary traces to the nonsemisimple setting; for strongly
rational orbifolds the corresponding statement holds with
ordinary traces \cite{HuangGenusOne}.

\paragraph{Proof strategy and organization.}
The main technical ingredient in this paper is a uniform derivative
bound in a relative $C_2$-quotient. Write $B=V^G$ and set
\[
 Q=R_B(V)=V/\Span_{\C}\{b_{-2}v\mid b\in B,\ v\in V\},
 \qquad q:V\longrightarrow Q.
\]
Choose a finite-dimensional graded $G$-stable space $U$ of strong
generators of $V$, and write $r=|G|$. By Noether's degree bound,
applied to an associated graded algebra, the construction of
invariant generators reduces to normally ordered products of derivatives
of at most $r$ elements of $U$. In other words, the degree, or number
of factors in each product, is at most $|G|$. The orders of differentiation, however,
may be arbitrarily large. Translation covariance gives
$(T^ju)_{-1}=j!\,u_{-1-j}$, so we seek a constant $H=H(U,r)$ such
that, for $k\le r$, homogeneous $u_i\in U$, and $j_i\ge0$,
\[
 q\bigl((u_1)_{-1-j_1}\cdots(u_k)_{-1-j_k}\ket{0}\bigr)=0
 \quad\text{whenever}\quad
 \sum_{i=1}^k j_i>H-\sum_{i=1}^k\wt u_i.
\]
The class on the left is the constant term in $z_1,\ldots,z_k$ of
\[
 \frac{1}{j_1!\cdots j_k!}
 \partial_{z_1}^{j_1}\cdots\partial_{z_k}^{j_k}
 q\!\left(Y(u_1,z_1)\cdots Y(u_k,z_k)\ket{0}\right),
\]
where $q$ acts coefficientwise. Thus the desired bound amounts to showing that derivatives of the
projected product have zero constant term in $Q$ for sufficiently
large total derivative order, uniformly in the $u_i$ and $k\le r$. Proposition~\ref{prop:cutoff}
gives such a bound for any fixed finite-dimensional graded
$G$-stable subspace of $V$ and any prescribed maximum number of insertions.

To that end, consider $n$-point functions
$Y(a_1,z_1)\cdots Y(a_n,z_n)\ket{0}$, where
$\mathbf{z}=(z_1,\ldots,z_n)$ and the $a_i$ lie in fixed
homogeneous $G$-stable input spaces $E_i\subset V_{h_i}$.
Multiplication by suitable locality factors gives a formal
power series $H_V(\mathbf{z})$; write $H_Q=q\circ H_V$ for its
image in $Q$. We first show that, after taking invariants in a chosen
subset of the inputs, sufficiently high derivatives of
$H_Q$ vanish when those inputs are translated together.
In particular, this gives a system of differential
equations for $H_Q$. However, the required orders
depend on the subset and the weights of the chosen
states, so these equations do not yet provide the desired uniform bound. More precisely, for a nonempty subset $I\subseteq[n]$, we
introduce collision coordinates $z_i=t+x_i$ for $i\in I$,
with one $x_i=0$. Locality and weak associativity
regroup the corresponding fields into an iterate in
the common coordinate $t$. Taking invariants in the
tensor factors indexed by $I$ makes the coefficients
of this iterate lie in $B$. The relative $C_2$-relations
then give the partial-translation equations
\[
 \partial_I^{d_I}(H_Q\circ P_I)=0,\qquad
 \partial_I=\sum_{i\in I}\partial_{z_i},\qquad
 d_I=1+\sum_{\substack{i\in I\\j\notin I}}(h_i+h_j),
\]
where $P_I$ projects onto invariants in the indicated
factors (Section~\ref{sec:fields}).

We encode these equations algebraically by replacing the commuting
differential operators $\partial_{z_i}$ with polynomial variables
$p_i$. Let $E=\bigotimes_iE_i$, $R=\C[p_1,\ldots,p_n]$,
and $p_I=\sum_{i\in I}p_i$. Define the universal module
\[
 M=\frac{R\otimes E}{
  \displaystyle\sum_{\varnothing\ne I\subseteq[n]}
    R\,p_I^{d_I}P_I E}.
\]
The coefficient map
\[
 R\otimes E\longrightarrow Q,\qquad
 \mathbf{p}^{\boldsymbol{\alpha}}\otimes e
 \longmapsto
 (\partial_{\mathbf{z}}^{\boldsymbol{\alpha}}H_Q)(\mathbf{0})(e)
\]
factors through $M$ by the partial-translation equations.

We now append specially chosen auxiliary fields, obtaining an
enlarged input space $\widetilde E$ and regularized field product
$\widetilde H_Q$. The auxiliary fields have two properties.
First, their representation-theoretic structure ensures that
the submodule $\mathcal N$ of the corresponding enlarged universal
module, generated by the prescribed tensors, satisfies
\[
 \Supp\mathcal N\subseteq\{\mathbf{0}\}.
\]
Since $\mathcal N$ is finitely generated, some power of the
maximal ideal at the origin annihilates it.
This gives a uniform bound on the total derivative order of
$\widetilde H_Q$ on the prescribed tensors. Second, the auxiliary fields occur in pairs
whose leading operator product coefficient is the vacuum.
Successively colliding these pairs and extracting their vacuum
coefficients transfers the bound to the original $H_Q$, proving
Proposition~\ref{prop:cutoff} (Section~\ref{sec:cutoff}).

Finally, the uniform derivative bound and Noether's degree bound
give finite generation of $R_B$, hence strong finite generation
of $B$ by Lemma~\ref{lem:stronggen-lift}
(Section~\ref{sec:generators}). The fixed-weight comparison with a decreasing filtration
then allows us to apply this bound to the kernel of
$R_B\to(R_V)^G$. Together with a differential
invariant-theoretic nilpotence argument, this proves
that the kernel is nilpotent (Section~\ref{sec:geometry}). Section~\ref{sec:consequences}
develops the applications, and the appendix supplies the
formal expansion and vacuum-coefficient arguments.

\paragraph{Acknowledgements} 
The author thanks Yi-Zhi Huang for his feedback and
references concerning twisted modules and genus-one
correlation functions, and Andrew Linshaw for drawing
attention to previous work on the vertex-algebra
Hilbert problem.

\paragraph{AI disclosure.}
Generative AI tools, including OpenAI's ChatGPT, were used
extensively under the author's direction to explore proof
strategies, develop mathematical arguments and revise exposition.
AI-assisted exploration contributed in particular to the
differential-equation and auxiliary-field arguments for the
uniform derivative bound (Sections~\ref{sec:fields}
and~\ref{sec:cutoff}), and to aspects of the strong-generation
and nilpotent-kernel proofs (Sections~\ref{sec:generators}
and~\ref{sec:geometry}). The author formulated and directed
the research, critically revised and independently checked the
arguments. 
The author takes full responsibility for the mathematical
claims, proofs, references, and final text.

\section{The relative \texorpdfstring{$C_2$}{C2} quotient}\label{sec:relative}
We begin by fixing notation and establishing several
preliminary results that will be used throughout the paper.

Set $B=V^G$ and replace $G$ by its image in $\operatorname{Aut}(V)$.
Write $r=|G|$ and $P_G=r^{-1}\sum_{g\in G}g$ for the projection onto
invariants. We use $Y(a,z)=\sum_{m\in\mathbb Z}a_mz^{-m-1}$, write
$T$ for translation and $\wt a$ for the degree of a homogeneous
state, and associate normally ordered products to the right.

\subsection{Relative relations and the comparison map}
For a vertex algebra $W$ and a $W$-module $M$, write
$R_W(M)=M/\Span\{b_{-2}v:b\in W,\ v\in M\}$. We apply this
construction to $V$ as a $B$-module, writing
\begin{equation}\label{eq:relative}
 D_B=\Span\{b_{-2}v:b\in B,\ v\in V\},\qquad
 Q=R_B(V)=V/D_B,\qquad q:V\to Q.
\end{equation}
The quotient is a graded $G$-module. Two properties account for its
use here. First, translation gives
$(T^{m-2}b)_{-2}=(m-1)!\,b_{-m}$ for $m\ge2$, so
\begin{equation}\label{eq:negative}
 q(b_{-m}v)=0 \qquad(b\in B,\ v\in V,\ m\ge2).
\end{equation}
Second, if $\sum_i b^i_{-2}v_i\in D_B$ is invariant, applying $P_G$
replaces each $v_i$ by $P_Gv_i$. It follows that
\begin{equation}\label{eq:exact}
 D_B\cap B=C_2(B),\qquad Q^G\simeq R_B.
\end{equation}
Indeed, every invariant class has an invariant lift by exactness.
Since $D_B\subseteq C_2(V)$, taking invariants in $Q\to R_V$ now gives
\begin{equation}\label{eq:comparison}
 0\longrightarrow K_G\longrightarrow R_B
 \xrightarrow{\ p_G\ }(R_V)^G\longrightarrow0,\qquad
 K_G=\frac{B\cap C_2(V)}{C_2(B)}.
\end{equation}
The product $[a_{-1}b]$ and bracket $[a_0b]$ are preserved by
inclusion, so $p_G$ is Poisson. We use $Q$ only as a vector space
with grading and $G$-action. All fields act in $V$ before applying
$q$.

\subsection{Weight and total derivative order}\label{subsec:strong}
A subspace $U\subset V$ strongly generates $V$ if the words
\begin{equation}\label{eq:words}
 W=u^1_{-n_1}\cdots u^\ell_{-n_\ell}\ket{0},
 \qquad u^i\in U,\quad n_i\ge1,
\end{equation}
span $V$, including the empty word. We say that $V$ is strongly
finitely generated if $U$ can be chosen finite dimensional. For
homogeneous $u^i$ of positive weight,
\begin{equation}\label{eq:wordweight}
 \wt W=\sum_i\wt u^i+\sum_i(n_i-1).
\end{equation}
The identity $(T^{n-1}u)_{-1}=(n-1)!\,u_{-n}$ identifies $n_i-1$
with the derivative order of the $i$th field. Their sum
$\sum_i(n_i-1)$ is the \emph{total derivative order} of this word.

\begin{lemma}[Strong-generation lifting criterion]\label{lem:stronggen-lift}
Let $W$ be a vertex algebra with a connected nonnegative integer
grading, and let $S\subseteq W_{>0}$ consist of homogeneous states.
Then $S$ strongly generates $W$ if and only if its images generate
$R_W=W/C_2(W)$ as a unital commutative algebra. In particular,
$W$ is strongly finitely generated if and only if $R_W$ is
finitely generated as a commutative algebra.
\end{lemma}
\begin{proof}
If $S$ strongly generates $W$, the spanning words
\eqref{eq:words} (with $W$ in place of $V$) map to products
of the images of $S$ in $R_W$, since modes $-n$ for $n\ge2$
vanish in $R_W$.

Conversely, put
\[
 C_1(W)=TW+\Span_{\C}\{a_{-1}b:a,b\in W_{>0}\}.
\]
We have $C_2(W)\subseteq C_1(W)$: for positive-weight $b$,
$a_{-2}b=(Ta)_{-1}b$, while $a_{-2}\ket{0}=Ta$.
If the images of $S$ generate $R_W$, every positive-weight
state is congruent modulo $C_2(W)$ to a polynomial in $S$.
Terms involving at least two positive-weight factors lie in
$C_1(W)$, so $S$ spans $W_{>0}/C_1(W)$. Li's
strong-generation criterion \cite[Theorem~4.10]{Li} now
shows that $S$ strongly generates $W$.

Finally, a finitely generated graded algebra $R_W$ admits a
finite homogeneous generating set, whose positive-weight
elements have homogeneous lifts to $W$. This proves the last
assertion.
\end{proof}

At fixed input weights, \eqref{eq:wordweight} turns an output-weight
bound into a bound on total derivative order. To obtain it uniformly
we keep the words as coefficients of field products and use
\eqref{eq:negative} before extracting those coefficients.

\section{Differential equations for partial translations}\label{sec:fields}

Before proving the uniform derivative bound, we establish
a system of differential equations for regularized
products of vertex operators modulo $Q$.
More precisely, after taking invariants in selected
input factors, sufficiently high derivatives in their
common translation direction vanish.
We prove this by using locality and weak associativity to
regroup the corresponding fields, and then applying
the defining relations of $Q$.
\subsection{Regularized \texorpdfstring{$n$}{n}-point functions}
Fix finite-dimensional $G$-stable subspaces $E_i\subset V_{h_i}$,
$1\le i\le n$. Set $E=\bigotimes_i E_i$ and
$\mathbf{z}=(z_1,\ldots,z_n)$. For a pure tensor
$e=a_1\otimes\cdots\otimes a_n$, define
\begin{align}
 \Delta(\mathbf{z})&=\prod_{i<j}(z_i-z_j)^{h_i+h_j},
       \label{eq:Delta}\\
 H_V(\mathbf{z})(e)&=\Delta(\mathbf{z})
       Y(a_1,z_1)\cdots Y(a_n,z_n)\ket{0},\qquad H_Q=q\circ H_V.
       \label{eq:HV}
\end{align}
Before regularization the product is expanded in
$V((z_1))\cdots((z_n))$ or equivalently $|z_1|>\cdots>|z_n|$. Locality clears its poles along the
diagonals, and creation accounts for regularity in each remaining
coordinate \cite{FHL,Frenkel}. The weight formula below converts a
polynomial degree bound for the regularized product into the output-weight
cutoff used in Proposition~\ref{prop:cutoff}. The proof is in
Appendix~\ref{app:regularity}.

\begin{lemma}[Regularity]\label{lem:regularity}
We have $H_V\in\operatorname{Hom}_{\C}(E,V[[z_1,\ldots,z_n]])$.
A coefficient of total degree $d$ has weight
\begin{equation}\label{eq:degreeweight}
 \sum_i h_i+d-\deg\Delta.
\end{equation}
In particular, projection to $V_w$ gives a homogeneous polynomial in $z_1,\ldots,z_n$.
\end{lemma}

Let $I=\{i_1<\cdots<i_m\}$ be a nonempty subset of
$[n]=\{1,\ldots,n\}$, with complement $I^c$. Translate the
coordinates in $I$ together by writing
\begin{equation}\label{eq:coordinates}
 z_i=t+x_i\quad(i\in I),\qquad x_{i_m}=0.
\end{equation}
Here $t=z_{i_m}$ and the $x_i$ record displacements from it.
Factor $\Delta=\Delta_I\Delta_{I^c}\Delta_{\rm cross}$, retaining
the orientation of each difference. The first two factors contain
the differences internal to $I$ and $I^c$, respectively. Only the
cross factor depends on $t$, and its $t$-degree is
\[
 D_I=\sum_{i\in I,\,j\notin I}(h_i+h_j).
\]
For $I^c=\{l_1<\cdots<l_{n-m}\}$, write the two internal products as
\begin{align*}
 C_I(\mathbf{x})&=\Delta_I(\mathbf{x})
    Y(a_{i_1},x_{i_1})\cdots Y(a_{i_{m-1}},x_{i_{m-1}})a_{i_m}
    =\sum_{\boldsymbol{\alpha}\ge\mathbf{0}}c_{\boldsymbol{\alpha}}
        \mathbf{x}^{\boldsymbol{\alpha}},\\
 W_{I^c}(\mathbf{z}_{I^c})&=\Delta_{I^c}(\mathbf{z}_{I^c})
    Y(a_{l_1},z_{l_1})\cdots Y(a_{l_{n-m}},z_{l_{n-m}})\ket{0}
    =\sum_{\boldsymbol{\gamma}\ge\mathbf{0}}v_{\boldsymbol{\gamma}}
        \mathbf{z}_{I^c}^{\boldsymbol{\gamma}}.
\end{align*}
The first is the regularized vacuum product for $I$ at $x_{i_m}=0$,
so both are power series by Lemma~\ref{lem:regularity}. For a
singleton $I$, $C_I=a_{i_1}$, and for $I^c=\varnothing$ we put
$W_{I^c}=\ket{0}$. We use the usual multi-index conventions,
in particular $|\boldsymbol{\alpha}|=\sum_i\alpha_i$ and
$\boldsymbol{\alpha}!=\prod_i\alpha_i!$, and compare multi-indices
componentwise. Brackets
$[\mathbf{x}^{\boldsymbol{\alpha}}]$ extract the indicated coefficient
and leave the remaining variables formal.

The following form of formal associativity keeps the cross factor
explicit. In the two-field case the product and iterate are different
expansions of one expression \cite[Proposition~2.6]{Frenkel}.
Appendix~\ref{app:iterate} verifies the corresponding reexpansion
for the indicated subset of fields, coefficient by coefficient.

\begin{lemma}[Reexpansion at a common coordinate]\label{lem:iterate}
In $V((t))[[\mathbf{x},\mathbf{z}_{I^c}]]$ one has
\begin{equation}\label{eq:iterateformula}
 \left.H_V(\mathbf{z})(e)\right|_{z_i=t+x_i\ (i\in I)}
 =
 \Delta_{\rm cross}(t+\mathbf{x},\mathbf{z}_{I^c})
       Y(C_I(\mathbf{x}),t)\,W_{I^c}(\mathbf{z}_{I^c}).
\end{equation}
The left side is substitution in a power series. Before multiplication
by the cross factor, the coefficient of
$\mathbf{x}^{\boldsymbol{\alpha}}\mathbf{z}_{I^c}^{\boldsymbol{\gamma}}$
on the right is $Y(c_{\boldsymbol{\alpha}},t)v_{\boldsymbol{\gamma}}$.
Its lower Laurent bound in $t$ may depend on both multi-indices.
\end{lemma}

\subsection{Translation after taking invariants}
Let $P_I$ take invariants for the diagonal action on the factors
of $E$ indexed by $I$, leaving the other factors fixed. When
$H_V$ is precomposed with $P_I$, equivariance replaces every
$c_{\boldsymbol{\alpha}}$ in \eqref{eq:iterateformula} by
$P_Gc_{\boldsymbol{\alpha}}\in B$. This is where the relative
relations enter.

\begin{lemma}[Partial translations]\label{lem:PDE}
For each nonempty $I\subseteq[n]$, put
$\partial_I=\sum_{i\in I}\partial_{z_i}$ and $d_I=D_I+1$. Then
\begin{equation}\label{eq:PDE}
 \partial_I^{d_I}(H_Q\circ P_I)=0 .
\end{equation}
\end{lemma}
\begin{proof}
For $b=P_Gc_{\boldsymbol{\alpha}}$, relation \eqref{eq:negative}
and lower truncation give
\[
 qY(b,t)v_{\boldsymbol{\gamma}}
   =\sum_{m\ge-1}q(b_mv_{\boldsymbol{\gamma}})t^{-m-1}
   \in Q[t^{-1}].
\]
Multiplication by $\Delta_{\rm cross}$ can therefore produce no
$t$-power above $D_I$, uniformly in the two multi-indices. On the
other side of \eqref{eq:iterateformula} regularity excludes negative
powers of $t$. The result is polynomial in $t$ of degree at most
$D_I$. With the displacements held fixed,
$\partial_t=\sum_{i\in I}\partial_{z_i}$, giving \eqref{eq:PDE}.
\end{proof}

Already for two fields there is a direction these equations can miss.
Take $G=\{\pm1\}$ and let both inputs lie in sign representations.
Then $P_{\{1\}}=P_{\{2\}}=0$, whereas $P_{\{1,2\}}$ is the identity
and $d_{\{1,2\}}=1$. The only equation is
$(\partial_{z_1}+\partial_{z_2})H_Q=0$, which admits power series
in $z_1-z_2$ of arbitrarily high degree. To bound the coefficients of
the actual field product, we obtain further equations from products
with additional insertions.

\section{Uniform bounds on derivatives}\label{sec:cutoff}

We now bound the weights of the full product, including the
dependence on coordinate differences left open by
Lemma~\ref{lem:PDE}. Here $V$ is simple and the input space may
be any finite-dimensional graded $G$-stable subspace.

\begin{proposition}[Uniform derivative bound]\label{prop:cutoff}
Let $X\subset V$ be a finite-dimensional graded $G$-stable subspace
and let $s\ge1$. There is an integer $H=H(X,s)\ge0$ such that,
for every $1\le k\le s$ and $a_1,\ldots,a_k\in X$, all
coefficients of
\[
 qY(a_1,z_1)\cdots Y(a_k,z_k)\ket{0}
\]
belong to $\bigoplus_{w=0}^H Q_w$.
\end{proposition}

To read this as a derivative bound, take homogeneous $a_i\in X$ of
weights $h_i$, with $1\le k\le s$, and put
$W=(a_1)_{-1-j_1}\cdots(a_k)_{-1-j_k}\ket{0}$ for $j_i\ge0$.
By \eqref{eq:wordweight}, $W$ has weight $\sum_i h_i+\sum_i j_i$. Hence
\[
 q(W)\ne0\quad\Longrightarrow\quad
 \sum_{i=1}^k j_i\le H(X,s)-\sum_{i=1}^k h_i .
\]
This bounds the sum of the derivative orders across all $k$ insertions
and hence each individual order $j_i$. For a word in strong generators,
this sum is the Li filtration index of its specified symbol in
Section~\ref{sec:geometry}.

The proof enlarges the product so that the equations from partial
translations apply to more subsets of positions. The additional
fields must do two things. At the level of $G$-modules their
prescribed tensor must lie in enough images of the projections
$P_I$ to force a degree bound. At the level of vertex operators it
must be possible to remove the extra insertions while retaining
that bound. Group differences supply the first property, and a
vacuum coefficient supplies the second.

\subsection{Choosing the auxiliary pairs}
The two identities below give different ways to place a tensor in
the images of subset projections. The first allows the subset to
lie anywhere among $r$ positions. In the second every subset
contains a prescribed position, at the cost of inserting the
vectors $(1-g)a_g$ in the other positions.

\begin{lemma}\label{lem:averaging}
For $r$ finite-dimensional $G$-modules $F_1,\ldots,F_r$,
\begin{equation}\label{eq:invariantsubsets}
 \bigotimes_{i=1}^r F_i
   =\sum_{\varnothing\ne I\subseteq[r]}\operatorname{im}P_I .
\end{equation}
For $G$-modules $F_0$ and $F_g$ indexed by $g\ne1$, $v\in F_0$,
and $a_g\in F_g$,
\begin{equation}\label{eq:difference}
 v\otimes\bigotimes_{g\ne1}(1-g)a_g
 \in\sum_{J\subseteq G\setminus\{1\}}\operatorname{im}P_{\{0\}\cup J}.
\end{equation}
Both statements may be tensored with arbitrary additional factors.
\end{lemma}
\begin{proof}
Give $A=\Sym(\bigoplus_{i=1}^r F_i)$ its multigrading by the
$r$ summands, and let $A_+$ be the ideal of elements of positive
total degree. For $x\in A_+$, the nonleading coefficients of the
orbit polynomial $\prod_{g\in G}(t-gx)$ belong to $(A_+)^G$.
Evaluating this polynomial at $x$ shows that
$x^r\in A(A_+)^G$. Polarizing over $\C$ gives
$A_+^r\subseteq A(A_+)^G$. In particular, this ideal contains
the whole multidegree-$(1,\ldots,1)$ component
$\bigotimes_iF_i$.

For a tensor $v_1\otimes\cdots\otimes v_r$, this gives an
expression as a sum of products $ab$, where $a\in A$ and
$b\in(A_+)^G$. Since $G$ preserves the multigrading, each
multihomogeneous component of $b$ is invariant. Taking the
multidegree-$(1,\ldots,1)$ component, we may therefore assume
$a$ and $b$ are multihomogeneous with multidegrees summing to
$(1,\ldots,1)$. Their supports partition $[r]$. Let $I$ be
the set of indices where $b$ has degree one; it is nonempty
because $b\in A_+$. Then
\[
 b\in\left(\bigotimes_{i\in I}F_i\right)^G,\qquad
 a\in\bigotimes_{i\notin I}F_i,
\]
so $ab=P_I(ab)\in\operatorname{im}P_I$. Summing over these
products proves \eqref{eq:invariantsubsets}.

For \eqref{eq:difference}, consider the identity
\begin{equation}\label{eq:masteridentity}
 \sum_{h\in G}hv\otimes\bigotimes_{g\ne1}(ha_g-ga_g)
 =v\otimes\bigotimes_{g\ne1}(a_g-ga_g).
\end{equation}
Indeed, when $h\ne1$, the factor indexed by $g=h$ in the
summand on the left is zero. Fix an order on $G\setminus\{1\}$.
For $J\subseteq G\setminus\{1\}$, define
\[
 b_g^J=
 \begin{cases}
 a_g,&g\in J,\\
 ga_g,&g\notin J.
 \end{cases}
\]
Expanding the left side of \eqref{eq:masteridentity} by choosing
$ha_g$ in the factors indexed by $J$ and $-ga_g$ in the
remaining factors gives
\[
 v\otimes\bigotimes_{g\ne1}(1-g)a_g
 =r\sum_{J\subseteq G\setminus\{1\}}
 (-1)^{r-1-|J|}P_{\{0\}\cup J}
 \left(v\otimes\bigotimes_{g\ne1}b_g^J\right).
\]
Here $P_{\{0\}\cup J}$ averages the diagonal $G$-action on
the indicated tensor factors, leaving all other factors fixed.
The factor $r$ comes from the normalization of this projection.
The sum includes $J=\varnothing$, whose term uses $P_{\{0\}}$;
the tensor factors retain their original order. Each summand
therefore belongs to the corresponding image, proving
\eqref{eq:difference}.
\end{proof}

To realize the group differences as auxiliary fields, choose a
homogeneous $a_g$ with $u_g=(1-g)a_g\ne0$ for each $g\ne1$.
Faithfulness gives such an $a_g$, and its weight is positive since
$G$ fixes $V_0$. We need a partner for $u_g$ whose leading product
with it is the vacuum. This is the sole use of simplicity in the
proof of the uniform derivative bound. Appendix~\ref{app:vacuum}
proves the following statement using the vacuum projection.

\begin{lemma}[Vacuum coefficients]\label{lem:detect}
For every nonzero $u\in V_h$, $h>0$, there is a homogeneous
$w\in V_k$, $k>0$, such that
\[
 u_{h+k-1}w\ne0\quad\text{in }\C\ket{0}.
\]
\end{lemma}

Choose and rescale $w_g$ so that
\begin{equation}\label{eq:pairunit}
 (u_g)_{h_g+k_g-1}w_g=\ket{0},\qquad h_g=\wt u_g,\quad k_g=\wt w_g .
\end{equation}
All higher modes vanish by the grading. Put $U_g=\Span(Ga_g)$ and
$W_g=\Span(Gw_g)$. Given $k$ original input spaces $E_i\subset V_{h_i}$
as in Section~\ref{sec:fields}, write $E=\bigotimes_i E_i$ and set
\[
 \widetilde E=E\otimes\bigotimes_{g\ne1}(U_g\otimes W_g)^{\otimes(2r-1)},
 \qquad
 \xi=\bigotimes_{g\ne1}(u_g\otimes w_g)^{\otimes(2r-1)} .
\]
There are $N=(r-1)(2r-1)$ pairs, written $(u_t,w_t)$ in this order
for $1\le t\le N$. The multiplicity ensures that after any $2r-2$
positions have been specified, a copy of every $u_g$ is still
available elsewhere. When $r=1$ there are no pairs and $\xi=1$.

\subsection{Support of the prescribed tensors}

For the $n=k+2N$ positions of $\widetilde E$, introduce a variable
$p_i$ for differentiation in the $i$th field coordinate. The
differential equations \eqref{eq:PDE} are encoded by the module
\begin{equation}\label{eq:M}
 R=\C[p_1,\ldots,p_n],\qquad p_I=\sum_{i\in I}p_i,\qquad
 M=\frac{R\otimes\widetilde E}
 {\displaystyle\sum_{\varnothing\ne I\subseteq[n]}
                   R\,p_I^{d_I}P_I\widetilde E}.
\end{equation}
Let $\mathcal N\subseteq M$ be the submodule generated by the classes
of $e\otimes\xi$ for $e\in E$. It is finitely generated since $E$
is finite dimensional. The two-field example at the end of
Section~\ref{sec:fields} explains why the subset projections by
themselves can leave support away from the origin. It suffices to
exclude that support for $\mathcal N$, since these are the tensors
from which the auxiliary pairs can be removed.

\begin{lemma}[Support at the origin]\label{lem:support}
For some integer $q_0\ge1$,
\begin{equation}\label{eq:origin}
 (p_1,\ldots,p_n)^{q_0}\mathcal N=0 .
\end{equation}
\end{lemma}
\begin{proof}
For $\boldsymbol{\lambda}=(\lambda_1,\ldots,\lambda_n)\in\C^n$,
write $\mathfrak m_{\boldsymbol{\lambda}}
=(p_1-\lambda_1,\ldots,p_n-\lambda_n)$.
We first prove that
$\mathcal N_{\mathfrak m_{\boldsymbol{\lambda}}}=0$
for every $\boldsymbol{\lambda}\ne\mathbf{0}$, and then use
finite generation to obtain the required power of
$(p_1,\ldots,p_n)$.

Fix such a $\boldsymbol{\lambda}$ and let
$S=\{i:\lambda_i\ne0\}$. By \eqref{eq:M}, the relations in
$M$ include $p_I^{d_I}P_Iv=0$ for every
$v\in\widetilde E$ and every nonempty $I\subseteq[n]$.
If $p_I(\boldsymbol{\lambda})\ne0$, then $p_I$ is a unit in
the local ring $R_{\mathfrak m_{\boldsymbol{\lambda}}}$, so
\[
 P_Iv=0\quad\text{in }M_{\mathfrak m_{\boldsymbol{\lambda}}}
 \qquad(v\in\widetilde E).
\]
We apply this observation in two cases.

Suppose first that $|S|\ge2r-1$. Choose a real line through
the origin in $\C$ containing none of the nonzero
$\lambda_i$. One of the corresponding open half-planes contains
at least $r$ of these numbers; choose their indices
$J\subseteq S$ with $|J|=r$. There exists a real-linear
functional $\ell:\C\to\mathbb R$ with
$\ell(\lambda_i)>0$ for every $i\in J$. Hence, for every
nonempty $I\subseteq J$,
\[
 \ell\bigl(p_I(\boldsymbol{\lambda})\bigr)
 =\sum_{i\in I}\ell(\lambda_i)>0.
\]
Thus $p_I(\boldsymbol{\lambda})\ne0$, and the image of
$P_I\widetilde E$ vanishes in
$M_{\mathfrak m_{\boldsymbol{\lambda}}}$ for each such $I$.
Applying \eqref{eq:invariantsubsets} to the $r$ factors indexed
by $J$, and tensoring with the remaining factors, gives
\[
 \widetilde E
 =\sum_{\varnothing\ne I\subseteq J}P_I\widetilde E.
\]
Since these tensors generate $M$ over $R$, we obtain
$M_{\mathfrak m_{\boldsymbol{\lambda}}}=0$ and therefore
$\mathcal N_{\mathfrak m_{\boldsymbol{\lambda}}}=0$.

Now suppose $1\le|S|\le2r-2$ and choose $j\in S$.
In our tensor order, the original factors occupy positions
$1,\ldots,k$, and the $t$th auxiliary pair occupies
positions $k+2t-1$ and $k+2t$, with $u_t$ first and
$w_t$ second. For each $g\ne1$, let $A_g\subseteq[n]$
be the set of positions occupied by the prescribed copies
of $u_g$ in $\xi$. These sets are pairwise disjoint, and
each has $2r-1$ elements. As $|S|\le2r-2$, we can choose
\[
 i_g\in A_g\setminus S\qquad(g\ne1).
\]
The selected positions are distinct, and none equals $j$.

By linearity it suffices to take a pure tensor
$e=e_1\otimes\cdots\otimes e_k\in E$. For each
$J\subseteq G\setminus\{1\}$, set
\[
 I_J=\{j\}\cup\{i_g:g\in J\}.
\]
Let $\theta_J\in\widetilde E$ be obtained from $e\otimes\xi$
by replacing the selected copy of $u_g=(1-g)a_g$ at
position $i_g$ with $a_g$ if $g\in J$ and with $ga_g$
otherwise, leaving all other factors unchanged. These
replacements lie in the specified $G$-stable spaces $U_g$.
The explicit identity \eqref{eq:masteridentity}, as expanded
in the proof of Lemma~\ref{lem:averaging}, applies with
distinguished position $j$ and the positions $i_g$.
It gives
\[
 e\otimes\xi
 =r\sum_{J\subseteq G\setminus\{1\}}
 (-1)^{r-1-|J|}P_{I_J}\theta_J.
\]
The factor at $j$ may be one of the original inputs or an
auxiliary field, since the identity allows an arbitrary vector
in that position. All other factors retain their original
tensor order, and the case $J=\varnothing$ is included.

Since every $i_g$ lies outside $S$, we have
\[
 p_{I_J}(\boldsymbol{\lambda})
 =\lambda_j+\sum_{g\in J}\lambda_{i_g}
 =\lambda_j\ne0.
\]
Our initial localization observation therefore makes every
$P_{I_J}\theta_J$ vanish in
$M_{\mathfrak m_{\boldsymbol{\lambda}}}$.
Thus every generator $e\otimes\xi$ of $\mathcal N$
vanishes after localization, proving
$\mathcal N_{\mathfrak m_{\boldsymbol{\lambda}}}=0$
also in this case.

If $\mathcal N=0$, take $q_0=1$. Otherwise, as
$\mathcal N$ is finitely generated over $R$,
$\Supp\mathcal N=V(\Ann_R\mathcal N)$.
We have excluded all maximal ideals except
$\mathfrak m_{\mathbf{0}}$ from this closed set.
The Nullstellensatz gives
\[
 (p_1,\ldots,p_n)\subseteq\sqrt{\Ann_R\mathcal N}.
\]
Choose integers $\nu_i\ge1$ such that
$p_i^{\nu_i}\mathcal N=0$ and put
\[
 q_0=1+\sum_{i=1}^n(\nu_i-1).
\]
Any monomial not divisible by one of the $p_i^{\nu_i}$
has total degree at most $\sum_i(\nu_i-1)$.
Consequently, every monomial of degree at least $q_0$
annihilates $\mathcal N$, proving \eqref{eq:origin}.
\end{proof}

\subsection{Removing the auxiliary pairs}

Lemma~\ref{lem:support} bounds the total order of derivatives on
the prescribed inputs. The next step transfers the resulting weight
bound through a collision of each auxiliary pair.

\begin{proof}[Proof of Proposition~\ref{prop:cutoff}]
By multilinearity fix $k\le s$ and homogeneous components
$E_i=X\cap V_{h_i}$ for the original fields. Use the enlarged
product \eqref{eq:HV} on $\widetilde E$, with variables in tensor
order
\[
 \boldsymbol{\zeta}=(z_1,\ldots,z_k,x_1,y_1,\ldots,x_N,y_N).
\]
Write $\widetilde\Delta$ for its regularizing polynomial
\eqref{eq:Delta}, $\widetilde H_Q$ for its image in $Q$, and
$H_t=\wt u_t$, $K_t=\wt w_t$ for the weights in pair $t$.

\smallskip\noindent\emph{The weight bound with the pairs inserted.}
Taylor coefficients define a linear map
\begin{equation}\label{eq:coeffmap}
 \Phi:R\otimes\widetilde E\longrightarrow Q,\qquad
 \mathbf{p}^{\boldsymbol{\alpha}}\otimes v\longmapsto
 (\partial_{\boldsymbol{\zeta}}^{\boldsymbol{\alpha}}
       \widetilde H_Q)(\mathbf{0})(v).
\end{equation}
This factors through $M$. Indeed, differentiating \eqref{eq:PDE}
and taking its constant coefficient kills every multiple of
$p_I^{d_I}\otimes P_Iv$. Put
$F(\boldsymbol{\zeta})(e)=\widetilde H_Q(\boldsymbol{\zeta})(e\otimes\xi)$.
Its $\boldsymbol{\alpha}$-coefficient is
$\Phi(\mathbf{p}^{\boldsymbol{\alpha}}\otimes(e\otimes\xi))/
\boldsymbol{\alpha}!$, hence vanishes for $|\boldsymbol{\alpha}|\ge q_0$
by \eqref{eq:origin}. This bounds the total degree of $F$ uniformly
in $e$. If $\widetilde h=\sum_{i=1}^k h_i+\sum_{t=1}^N(H_t+K_t)$, the
degree--weight identity \eqref{eq:degreeweight} gives
\begin{equation}\label{eq:weightbound}
 F_w=0\quad\text{for }w>H,\qquad
 H=\max\{0,\widetilde h+q_0-1-\deg\widetilde\Delta\},
\end{equation}
where the subscript denotes projection to $Q_w$.

\smallskip\noindent\emph{The effect of one collision.}
Let $H_V^{(j)}$ be the regularized product of the original fields
and the remaining pairs $j+1,\ldots,N$, in their original order.
Put $H_Q^{(j)}=qH_V^{(j)}$, so $H_Q^{(0)}=F$ and $H_Q^{(N)}=H_Q$.
For $t=j+1$ write $x_t=y_t+d$. The internal factor of this pair in
the regularizer is $d^{H_t+K_t}$. Write $B_t^{(j)}(d)$ for the
product of its cross factors with all remaining fields. By
\eqref{eq:pairunit} and the grading,
\[
 d^{H_t+K_t}Y(u_t,d)w_t=\ket{0}+d\rho_t(d),\qquad
 \rho_t(d)\in V[[d]].
\]
Reexpanding this pair by Lemma~\ref{lem:iterate} gives
\[
 \left.H_Q^{(j)}\right|_{x_t=y_t+d}
 =B_t^{(j)}(d)\,
 q\!\left(Y(\ket{0}+d\rho_t(d),y_t) H_V^{(j+1)}\right).
\]
The expression on the right has only nonnegative powers of $d$.
Its constant term is $B_t^{(j)}(0)H_Q^{(j+1)}$, since the vacuum
field acts as the identity. On the left regularity makes the same
operation the specialization $x_t=y_t$.

\smallskip\noindent\emph{Returning to the original product.}
Iterating the one-pair identity yields
\begin{equation}\label{eq:collision}
 \left.F\right|_{x_1=y_1,\ldots,x_N=y_N}
   = C(\mathbf{z},\mathbf{y})\,H_Q(\mathbf{z}),
\end{equation}
where the accumulated cross factors are
\begin{equation}\label{eq:C}
 C(\mathbf{z},\mathbf{y})=
 \prod_{i=1}^k\prod_{t=1}^N(z_i-y_t)^{2h_i+H_t+K_t}
 \prod_{1\le t<v\le N}
 (y_t-y_v)^{2(H_t+K_t+H_v+K_v)} .
\end{equation}
For an original field and a pair there are two cross factors, and
between distinct pairs there are four. Their exponents add as in
\eqref{eq:C}, while the factors among the original fields remain
in $H_Q$. Only the coordinates \emph{within} a pair were identified.
The different $y_t$ are independent, so $C$ is nonzero, with $C=1$
if $N=0$.

Fix $w>H$. By \eqref{eq:weightbound} the left side of
\eqref{eq:collision} vanishes in weight $w$. For each $e\in E$,
$(H_Q)_w(e)$ is a $Q_w$-valued polynomial by
Lemma~\ref{lem:regularity}, so multiplication by the nonzero scalar
polynomial $C$ is injective. Hence $(H_Q)_w(e)=0$.
Multiplication by $\Delta^{-1}$, expanded in the original ordered
Laurent space, commutes with the projection to $Q_w$. The same
vanishing therefore holds before regularization. Finally $X$ has
only finitely many weights and $k\le s$. Taking the maximum over
these finitely many choices gives $H(X,s)$.
\end{proof}

\begin{remark}
The choice in \eqref{eq:pairunit} used only one vacuum coefficient
for each nonidentity element of the effective group. Thus in
Proposition~\ref{prop:cutoff} and Theorems~\ref{thm:av}
and~\ref{thm:sfg} simplicity may be replaced by the existence,
for each such $g$, of homogeneous $a_g,w_g$
such that $u_g=(1-g)a_g\ne0$ and
\[
 (u_g)_{\wt u_g+\wt w_g-1}w_g\in\C^\times\ket{0}.
\]
\end{remark}

\section{Descent of Strong Finite Generation}\label{sec:generators}

Our goal in this section is to prove that strong finite generation descends to finite-group orbifolds. More precisely, if $V$ is strongly finitely generated and satisfies the standing hypotheses of Section~\ref{sec:intro}, then $B=V^G$ is also strongly finitely generated. By Lemma~\ref{lem:stronggen-lift}, it suffices to prove that
\[
R_B=B/C_2(B)
\]
is finitely generated as a commutative algebra. The main ingredient is the uniform derivative bound of Proposition 4.1. To apply it, however, we first need a uniform bound on the number of factors appearing in the normally ordered words required to generate $R_B$.

We obtain such a bound using the increasing generating-weight
filtration $F_\bullet V$, defined by the sum of the weights of the
strong generators appearing in a word. More precisely, for a finite-dimensional, graded, $G$-stable strong generating space $U$, we consider invariant parts of words of the form
\[
\beta=P_G\left((T^{j_1}u_1)_{-1}\cdots
(T^{j_\ell}u_\ell)_{-1}\ket{0}\right),
\qquad u_i\in U,\quad j_i\geq 0.
\]
Here $\ell$ is the number of factors, while the filtration degree $p$ and conformal weight $w$ are related by
\[
p=\sum_{i=1}^{\ell}\wt u_i,
\qquad
w=p+\sum_{i=1}^{\ell}j_i.
\]

Passing to the associated graded takes us to a commutative differential algebra, where we can apply classical invariant theory. In particular, Noether's bound states that the invariant ring of a polynomial algebra under a finite group $G$ is generated by invariants of polynomial degree at most $|G|$. This allows us to choose generating words with
\[
\ell\leq |G|,
\]
independently of their derivative orders. We can therefore apply Proposition 4.1 with $s=|G|$, obtaining a uniform bound $H=H(U,|G|)$ such that
\[
q(\beta)\neq 0
\quad\Longrightarrow\quad
w\leq H,\qquad
\sum_{i=1}^{\ell}j_i=w-p\leq H-p.
\]
Thus both the number of factors and their total derivative order are bounded. Since $U$ is finite-dimensional, only finitely many such generating words remain. Their classes generate $R_B$, so Lemma~\ref{lem:stronggen-lift}
gives strong finite generation of $V^G$.
\subsection{The generating-weight filtration}
Take $V$ simple and strongly finitely generated. Choose a finite-dimensional
graded $G$-stable strong generating space $U\subset V_{>0}$ by
taking the $G$-span of homogeneous generators. Write $U_h=U\cap V_h$,
set $F_pV=0$ for $p<0$, and let
\begin{equation}\label{eq:goodfiltration}
 F_pV=\Span\left\{
 u^1_{-n_1}\cdots u^\ell_{-n_\ell}\ket{0}:
 n_i\ge1,\ u^i\in U\text{ homogeneous},\
 \sum_i\wt u^i\le p\right\}\qquad(p\ge0).
\end{equation}
The empty word gives $F_0V=\C\ket{0}$. In this filtration $T^ju$
has the same filtration bound as $u$ and weight $\wt u+j$.
Li's good-filtration theorem \cite[Theorem~4.11]{Li} gives
\[
 (F_pV)_mF_qV\subseteq
 \begin{cases}
 F_{p+q}V,&m<0,\\
 F_{p+q-1}V,&m\ge0 .
 \end{cases}
\]
Thus $A_F=\gr_F V$ is a commutative differential algebra. It is
generated by $D^j\bar u$ for homogeneous $u\in U$ and $j\ge0$,
where $\bar u$ is the symbol in degree $\wt u$ and $D$ is induced
by $T$.

Restrict $F$ to $B$ and give $R_B$ the quotient filtration.
In weight $w$, $(F_pR_B)_w$ is the image of $F_pV\cap B_w$
in $(R_B)_w$. An invariant symbol in $A_F$ has an invariant
lift by applying $P_G$. Taking the quotient with the induced
filtration on $C_2(B)$ then gives
\begin{equation}\label{eq:filteredRB}
 \gr_F B=A_F^G,\qquad C_F:=\gr_F R_B=A_F^G/\gr_F C_2(B).
\end{equation}
We write $(C_F)_{p,w}=(F_pR_B)_w/(F_{p-1}R_B)_w$ for the
component of filtration degree $p$ and original weight $w$.
For $u\in U_h$, the specified symbol $D^j\bar u$ in $A_F$ has
bidegree $(h,h+j)$ whenever it is nonzero.

\subsection{Applying Noether's bound}
The filtration in \eqref{eq:filteredRB} serves two purposes.
Generators of $C_F$ lift to generators of $R_B$. Their bidegrees
will also measure how many derivatives can survive in the
comparison kernel.

\begin{proposition}\label{prop:filteredgenerators}
Let $H_U=H(U,r)$ from Proposition~\ref{prop:cutoff}. The algebra
$C_F$ is generated by the symbols of invariant parts of normally
ordered products of derivatives of at most $r$ elements of $U$,
of total weight at most $H_U$. These symbols constitute a finite
bihomogeneous generating set when the inputs range over a fixed
homogeneous basis of $U$.
\end{proposition}
\begin{proof}
For each $j\ge0$, take a copy $U^{(j)}$ of the $G$-module $U$.
There is a $G$-equivariant surjection
\[
 S=\Sym\left(\bigoplus_{j\ge0}U^{(j)}\right)
 \longrightarrow A_F,\qquad u^{(j)}\longmapsto D^j\bar u .
\]
It preserves both degrees and induces a surjection $S^G\to A_F^G$.
Every polynomial belongs to some
$S_{\le J}=\Sym(\bigoplus_{j=0}^J U^{(j)})$. Noether's theorem
\cite{Noether} gives generators of $S_{\le J}^G$ of polynomial degree
at most $r$, independently of $J$. Taking the union over $J$ shows
that $S^G$ is generated by $P_G\mathsf{m}$ for basis monomials
$\mathsf{m}=u_1^{(j_1)}\cdots u_\ell^{(j_\ell)}$ with
$1\le\ell\le r$. For $h_i=\wt u_i$, write the number of factors,
the sum of the generating weights, the total derivative order, and
the original weight as
\begin{equation}\label{eq:worddegrees}
 \deg_{\rm pol}\mathsf{m}=\ell,\qquad p=\sum_i h_i,\qquad
 j=\sum_i j_i,\qquad w=p+j .
\end{equation}
Noether bounds $\ell$. Even with $\ell$ and $p$ fixed, the
derivative orders can make $w$ arbitrarily large.

Fix an order on the basis symbols. Lift $P_G\mathsf{m}$ to
\begin{equation}\label{eq:invariantlift}
 \beta_{\mathsf{m}}=P_G\!\left(
    (T^{j_1}u_1)_{-1}\cdots(T^{j_\ell}u_\ell)_{-1}\ket{0}\right)
     \in F_pB\cap B_w .
\end{equation}
Its degree-$p$ symbol is the image of $P_G\mathsf{m}$ by the
good-filtration product rule. This statement allows that symbol to
be zero. The specified symbols generate $A_F^G$ and hence their images
generate $C_F$.

The lift in \eqref{eq:invariantlift} equals $\prod_i j_i!$ times
the invariant part of
\[
 [z_1^{j_1}\cdots z_\ell^{j_\ell}]
     Y(u_1,z_1)\cdots Y(u_\ell,z_\ell)\ket{0}.
\]
There are at most $r$ insertions, whatever the $j_i$. Hence
Proposition~\ref{prop:cutoff} and \eqref{eq:exact} give
\begin{equation}\label{eq:highliftzero}
 w>H_U\quad\Longrightarrow\quad q(\beta_{\mathsf{m}})=0
 \quad\Longrightarrow\quad \beta_{\mathsf{m}}\in C_2(B).
\end{equation}
Its specified symbol then vanishes in $C_F$. For $w\le H_U$, only
finitely many basis words and derivative orders occur. Their images
therefore generate the particular associated graded algebra in
\eqref{eq:filteredRB}.
\end{proof}
We may now show finite strong generation descends to $V^G$.
\begin{restatedsfg}[Finite strong generation]
Let $V$ be a simple, strongly finitely generated vertex algebra
with a connected nonnegative integer grading, and let $G$ be a
finite group of grading-preserving automorphisms. Then $V^G$
is strongly finitely generated.
\end{restatedsfg}
\begin{proof}
Lift the finite generating symbols of
Proposition~\ref{prop:filteredgenerators} to $R_B$ using
\eqref{eq:invariantlift}. They generate $R_B$ by induction on
filtration degree. Indeed, express the symbol of a class in
$(F_pR_B)_w$ as a bihomogeneous polynomial in these generators.
Subtract the same polynomial in the lifts to obtain a class in
$(F_{p-1}R_B)_w$. The induction starts with the constants in
$F_0R_B$. Lemma~\ref{lem:stronggen-lift} now gives finite strong
generation of $B$.
\end{proof}

\section{Associated varieties of finite-group orbifolds}\label{sec:geometry}

We now prove Theorem~\ref{thm:av}. Write $B=V^G$ and
$K_G=\ker(p_G:R_B\to(R_V)^G)$. Our goal is to show that
$K_G$ is nilpotent.

The proof uses two filtrations on $V$, which serve different
purposes. The increasing generating-weight filtration $F_\bullet V$
from Section~\ref{sec:generators} was used to construct a finite
set of generators of $R_B$. We now use the canonical decreasing
filtration $E^\bullet V$, which records derivative order at fixed
weight and satisfies $E^1V=C_2(V)$.
The two filtrations are related in each homogeneous weight by
Arakawa's comparison \cite[Proposition~2.6.1]{ArakawaC2}:
\begin{equation}\label{eq:Li-comparison}
 E^jV\cap V_w=F_{w-j}V\cap V_w .
\end{equation}
Thus the bounds on generating symbols obtained in Section~\ref{sec:generators}
also give bounds in the decreasing filtration.

To detect $K_G$, we use the filtration on $R_B$ induced by
$B\cap E^jV$, rather than the intrinsic decreasing filtration
of $B$. We first establish nilpotence in the associated graded
algebra and then use the uniform derivative bound to lift this
nilpotence to $R_B$.

\subsection{The decreasing filtration and the fixed-weight comparison}

The canonical decreasing filtration, introduced by Li
\cite{LiAbelianizing}, is defined by
\begin{equation}\label{eq:Li-definition}
 E^jV=\Span_{\C}\left\{
 a^1_{-1-m_1}\cdots a^s_{-1-m_s}b:
 a^i,b\in V,\ m_i\ge0,\ \sum_i m_i\ge j
 \right\}\qquad(j\ge0).
\end{equation}
Unlike the generating-weight filtration, this definition allows
arbitrary states $a^i,b\in V$, not only the chosen strong generators.
We have $E^0V=V$ and $E^1V=C_2(V)$. Translation takes $E^jV$
into $E^{j+1}V$, and $(E^iV)_{-1}(E^jV)\subseteq E^{i+j}V$
\cite[Theorem~2.12 and Corollary~3.6]{LiAbelianizing}.

The fixed-weight comparison \eqref{eq:Li-comparison} can also be
seen directly from the defining words. By \eqref{eq:worddegrees},
a defining $U$-word of weight $w$ with generating-weight sum at most
$w-j$ has at least $j$ derivatives. Conversely, in a homogeneous
word in \eqref{eq:Li-definition} of weight $w$, the weights of
$a^1,\ldots,a^s,b$ sum to $w-\sum_i m_i\le w-j$. Each input of
weight $h$ lies in $F_hV$ by \eqref{eq:wordweight}, so the
negative-mode product rule puts the word in $F_{w-j}V$. In particular
$E^jV\cap V_w=0$ for $j>w$.

Set $A=\gr_E V$, with
$A_{j,w}=(E^jV\cap V_w)/(E^{j+1}V\cap V_w)$.
Write $A_j=\bigoplus_w A_{j,w}$ and $A_{>0}=\bigoplus_{j>0}A_j$.
The differential generation theorem
\cite[Corollary~4.3]{LiAbelianizing} gives
\begin{equation}\label{eq:derivativeideal}
 A_0=R_V,\qquad A_{>0}=A\,D(A),
\end{equation}
where $D$ is induced by $T$. Indeed, $D$ raises $E$-degree by one,
and $A$ is differentially generated by $A_0$. Thus every monomial
of positive $E$-degree contains a derivative.

\subsection{Nilpotence of the comparison kernel}
The following differential identity controls products of
positive-degree symbols in the associated graded algebra. The
uniform derivative bound controls the degree of each generating symbol.

\begin{lemma}\label{lem:differential}
Let $A$ be a commutative $\C$-algebra with derivation $D$.
If a group $G$ of order $r$ acts by algebra automorphisms
commuting with $D$, then
\begin{equation}\label{eq:differentialideal}
 (A\,D(A))^{2r-1}\subseteq A\,D(A^G).
\end{equation}
\end{lemma}
\begin{proof}
The ideal $I=A\,D(A^G)$ is $D$-stable since $Df\in A^G$ for $f\in A^G$. Every invariant in
$\overline A=A/I$ lifts to $A^G$, so the induced derivation kills
$\overline A^G$. For $a\in\overline A$, let
$P(t)=\prod_{g\in G}(t-ga)$. Its coefficients are invariant.
Differentiating $P(a)=0$ gives $P'(a)Da=0$. We claim that
\[
 P^{(j)}(a)(Da)^{2j-1}=0\qquad(1\le j\le r),
\]
where $P^{(j)}$ is the $j$th derivative in $t$. Applying $D$ to
the $j$th identity and multiplying by $Da$ gives
\[
 P^{(j+1)}(a)(Da)^{2j+1}
 +(2j-1)P^{(j)}(a)(Da)^{2j-1}D^2a=0 .
\]
The second term vanishes by induction. At $j=r$ we have
$P^{(r)}=r!$, so $(Da)^{2r-1}=0$ for every $a\in\overline A$.
Polarizing over $\C$ shows that every product of $2r-1$ derivatives
in $\overline A$ vanishes. This proves \eqref{eq:differentialideal}.
\end{proof}
We are now ready to prove the main result of the paper.
\begin{restatedav}[Associated varieties of finite-group orbifolds]
Let $V$ be a simple, strongly finitely generated vertex algebra with
a connected nonnegative integer grading, and let $G$ be a finite
group of grading-preserving automorphisms. The canonical Poisson
homomorphism $p_G:R_{V^G}\twoheadrightarrow(R_V)^G$ has nilpotent
kernel. In particular it induces an isomorphism of reduced Poisson
varieties
\[
 X_V/G\ \xrightarrow{\ \sim\ }\ X_{V^G}.
\]
\end{restatedav}
\begin{proof}
\emph{ Identify the kernel in the inherited filtration.}
Give $B=V^G$ the filtration $B\cap E^jV$ inherited from $V$, and
denote the resulting quotient filtration on $R_B$ by $\mathcal E$.
By $E^1V=C_2(V)$ and \eqref{eq:comparison},
\begin{equation}\label{eq:inherited}
 \mathcal E^jR_B=\frac{(B\cap E^jV)+C_2(B)}{C_2(B)},\qquad
 \mathcal E^1R_B=\frac{B\cap C_2(V)}{C_2(B)}=K_G .
\end{equation}
The first step of the intrinsic decreasing filtration of $B$ is
$C_2(B)$, whose image in $R_B$ is zero. The $(-1)$-product rule for $E$ gives
$\mathcal E^iR_B\,\mathcal E^jR_B\subseteq\mathcal E^{i+j}R_B$.

Put $C=\gr_{\mathcal E}R_B$. Restricting \eqref{eq:Li-comparison}
to $B_w$ and taking images in $(R_B)_w$ gives
$(\mathcal E^jR_B)_w=(F_{w-j}R_B)_w$. Passing to consecutive
quotients therefore identifies the two gradings,
\begin{equation}\label{eq:regraded-quotient}
 C_{j,w}
 =\frac{(\mathcal E^jR_B)_w}{(\mathcal E^{j+1}R_B)_w}
 =\frac{(F_{w-j}R_B)_w}{(F_{w-j-1}R_B)_w}
 =(C_F)_{w-j,w}.
\end{equation}
Both products are induced by multiplication in $R_B$, so the
generating symbols of Proposition~\ref{prop:filteredgenerators}
generate $C$ after this change of indices. As for $A$, write
$C_j=\bigoplus_w C_{j,w}$ and $C_{>0}=\bigoplus_{j>0}C_j$.

\smallskip\noindent
\emph{Nilpotence in the associated graded algebra.}
Exactness of invariants identifies $(A_j)^G$ with
$(B\cap E^jV)/(B\cap E^{j+1}V)$ in each $E$-degree. By the
definition of $\mathcal E$ this space maps onto $C_j$. Together
these maps give a bigraded algebra surjection $\pi:A^G\to C$.
For a bihomogeneous invariant symbol $f$, choose an invariant lift
$b\in B$. In the next $E$-degree, $Df$ is represented by
$Tb=b_{-2}\ket{0}\in C_2(B)$, so $\pi(Df)=0$.
Taking invariant parts of the coefficients in a
sum $\sum_i a_iDf_i$ with $f_i\in A^G$ gives
\[
 (A\,D(A^G))^G=A^G D(A^G)\subseteq\ker\pi .
\]
Every class of positive $\mathcal E$-degree in $C$ lifts
in each bidegree
to $A^G\cap A_{>0}$. By \eqref{eq:derivativeideal} and
Lemma~\ref{lem:differential}, the product of $2r-1$ such lifts
belongs to $(A\,D(A^G))^G$. Hence
\begin{equation}\label{eq:gradednil}
 (C_{>0})^{2r-1}=0 .
\end{equation}

\smallskip\noindent
\emph{Lift back to $R_B$.}
Use the finite set of generating symbols from
Proposition~\ref{prop:filteredgenerators}, reindexed by
\eqref{eq:regraded-quotient}. Each has $\mathcal E$-degree $j=w-p\le w\le H_U$,
where $H_U=H(U,r)$. Let $d$ be the largest positive
$\mathcal E$-degree of the nonzero symbols in this set,
with $d=0$ if there is none.
By \eqref{eq:gradednil} a nonzero monomial contains at most $2r-2$ positive-degree
generators. Each contributes at most $d$, and generators of
degree zero contribute nothing. Consequently
\begin{equation}\label{eq:deficitbound}
 C_j=0\qquad(j>M),\qquad M=(2r-2)d .
\end{equation}

Fix an original weight $w$. By \eqref{eq:regraded-quotient}
and \eqref{eq:deficitbound}, the spaces $(\mathcal E^jR_B)_w$
and $(\mathcal E^{j+1}R_B)_w$ agree for every $j>M$. Choose
$J=\max\{M+1,w+1\}$, so $E^JV\cap V_w=0$ by
\eqref{eq:Li-comparison}. These equalities give
\[
 (\mathcal E^{M+1}R_B)_w=(\mathcal E^JR_B)_w=0 .
\]
This holds in every weight with the same $M$. By multiplicativity
and \eqref{eq:inherited}, the kernel in the original algebra satisfies
\begin{equation}\label{eq:kernelexponent}
 K_G^{\,1+(2r-2)d}=(\mathcal E^1R_B)^{M+1}
       \subseteq\mathcal E^{M+1}R_B=0 .
\end{equation}

Since $p_G$ is surjective with nilpotent kernel, it induces the first
isomorphism below. The second follows by taking invariants in
$R_V\to (R_V)_{\mathrm{red}}$, since the nilradical of $(R_V)^G$
is $(R_V)^G\cap\sqrt{(0)}$ in $R_V$. We obtain
\[
 (R_B)_{\mathrm{red}}\ \xrightarrow{\ \sim\ }\
 ((R_V)^G)_{\mathrm{red}}\ \xrightarrow{\ \sim\ }\
 ((R_V)_{\mathrm{red}})^G .
\]
In characteristic zero the nilradical is stable under derivations,
in particular under Hamiltonian derivations. Taking spectra therefore
gives the canonical Poisson isomorphism $X_V/G\simeq X_B$.
\end{proof}

\section{Consequences for finite-group orbifolds}\label{sec:consequences}

We conclude by deriving several consequences of the description
$X_{V^G}\simeq X_V/G$ and the descent of strong finite generation.
The first concerns the Poisson geometry of associated varieties;
the remaining applications concern the representation theory of
vertex operator algebras.

\subsection{Quasi-lisse orbifolds and ordinary characters}
\label{subsec:leaves}

Theorem~\ref{thm:av} identifies $X_{V^G}$ with the
finite Poisson quotient $X_V/G$. Since finite Poisson
quotients preserve the finiteness of symplectic leaves
\cite{ES2010}, we obtain the following consequence.

\begin{restatedquasilisse}[Quasi-lisse descent]
Let $V$ be a simple quasi-lisse vertex algebra with
a connected nonnegative integer grading, and let
$G$ be a finite group of grading-preserving automorphisms.
Then $V^G$ is quasi-lisse.
\end{restatedquasilisse}
\begin{proof}
By Theorem~\ref{thm:sfg}, the vertex algebra $V^G$ is strongly
finitely generated. Since $G$ acts on $X_V$ by Poisson
automorphisms, the finite quotient $X_V/G$ has finitely many
symplectic leaves whenever $X_V$ does \cite{ES2010}.
Theorem~\ref{thm:av} identifies this quotient with $X_{V^G}$.
\end{proof}

This has representation-theoretic consequences. If $V$ is a simple
quasi-lisse VOA of CFT type, then so is $B=V^G$. By the theorems
of Arakawa and Kawasetsu \cite{AK}, $B$ has only finitely many
simple ordinary modules, and the normalized characters of its
ordinary modules satisfy a common modular linear differential
equation. These conclusions do not require $B$ to be rational or
$C_2$-cofinite.

\subsection{$C_2$-cofiniteness and strong rationality}\label{subsec:rational}

The zero-dimensional case gives a stronger finiteness property.
For a strongly finitely generated vertex algebra $W$ with connected
nonnegative grading, Lemma~\ref{lem:stronggen-lift} shows that
$R_W$ is a finitely generated commutative algebra. Consequently,
$X_W$ is a point if and only if $R_W$ is finite dimensional.
Thus the condition that $X_W$ is a point is equivalent to
$C_2$-cofiniteness.

\begin{corollary}[$C_2$-cofiniteness]\label{cor:c2}
Let $V$ be a simple vertex algebra with a connected nonnegative
integer grading, and let $G$ be a finite group of grading-preserving
automorphisms. If $V$ is $C_2$-cofinite, then $V^G$ is
$C_2$-cofinite.
\end{corollary}
\begin{proof}
Since $R_V$ is finite dimensional, Lemma~\ref{lem:stronggen-lift}
shows that $V$ is strongly finitely generated. Its associated
variety is a point, so Theorem~\ref{thm:av} gives that
$X_{V^G}$ is a point. By Theorem~\ref{thm:sfg}, the orbifold is
strongly finitely generated. The preceding criterion therefore
implies that $R_{V^G}$ is finite dimensional.
\end{proof}

We can now complete the deduction of strong rationality.

\begin{restatedrational}[Descent of Strong Rationality]
Let $V$ be a strongly rational VOA and
$G$ a finite group of VOA automorphisms.
Then $V^G$ is strongly rational.
\end{restatedrational}
\begin{proof}
McRae proved that if $V$ is strongly rational and $G$ is finite,
then $V^G$ is strongly rational provided it is $C_2$-cofinite
\cite[Corollary~4.23]{McRae}. The latter condition follows from
Corollary~\ref{cor:c2}.
\end{proof}

In particular, when $V$ is strongly rational, the category of
ordinary $V^G$-modules is a modular tensor category
\cite{HuangRigidity}. For the tensor-category structure of twisted
sectors, the weaker hypothesis of $C_2$-cofiniteness already
suffices.

\subsection{Twisted sectors and holomorphic orbifolds}\label{subsec:twisted}

Let $V$ now be a simple $C_2$-cofinite VOA of CFT type and
$G\leq\operatorname{Aut}(V)$ finite. A grading-restricted
generalized $g$-twisted $V$-module restricts naturally to a
$V^G$-module, since every element of $V^G$ has trivial
$g$-monodromy. We restate the resulting categorical structure
from the introduction.

\begin{restatedtwisted}[Twisted sectors of a finite orbifold]
Let $V$ be a simple $C_2$-cofinite VOA of CFT type and
$G\leq\operatorname{Aut}(V)$ a finite group. Then
\[
 \operatorname{TwMod}_G(V)
 =\bigoplus_{g\in G}\operatorname{TwMod}_g(V)
\]
admits a natural $G$-crossed braided tensor category structure.
Its equivariantization is braided tensor equivalent to
$\operatorname{Mod}(V^G)$.
\end{restatedtwisted}
\begin{proof}
Corollary~\ref{cor:c2} makes $B=V^G$ $C_2$-cofinite, and $B$
inherits CFT type from $V$. Huang's cofiniteness theorem
\cite{HuangCofinite}, together with logarithmic tensor-category
theory \cite{HLZ}, supplies the requisite braided vertex tensor
category structure on $\operatorname{Mod}(B)$. McRae's
Theorems~4.15 and~4.17 \cite{McRae} then give the
$G$-crossed braided category of twisted modules and identify its
$G$-equivariantization with $\operatorname{Mod}(B)$.
\end{proof}

In the strongly rational case, related $G$-crossed tensor
structures and crossed modularity are studied by Xu
\cite{XuCrossed,XuCrossedModularity}.

\paragraph{Holomorphic orbifolds and reconstruction.}\label{subsec:holomorphic}
Theorem~\ref{thm:rational} removes the orbifold-rationality
assumption in the reconstruction theorem of Gannon and the
author \cite[Theorem~3.2]{GR}, building on Evans--Gannon
\cite{EG}. Consequently, for every finite group $G$ and
class $[\omega]\in H^3(G,\C^\times)$, there is a holomorphic
VOA $V$ with a $G$-action such that
\[
 \operatorname{TwMod}_G(V)\simeq\mathrm{Vec}_G^\omega,
 \qquad
 \operatorname{Mod}(V^G)\simeq\mathcal Z(\mathrm{Vec}_G^\omega),
\]
as tensor and braided tensor categories, respectively.
The reconstruction is that of \cite{EG,GR}; our result
supplies the previously conditional orbifold-rationality input.

\subsection{Genus-one modular invariance}\label{subsec:genusone}

A final consequence follows from $C_2$-cofiniteness without any
semisimplicity assumption. Huang's modular invariance theorem
for logarithmic intertwining operators \cite{HuangLog} applies to
the fixed-point VOA by Corollary~\ref{cor:c2}.

\begin{corollary}[Genus-one modular invariance]\label{cor:genusone}
Let $V$ be a simple $C_2$-cofinite VOA of CFT type, let
$G\leq\operatorname{Aut}(V)$ be finite, and set $B=V^G$.
Fix $n\geq 1$, grading-restricted generalized $B$-modules
$W_1,\ldots,W_n$, and vectors $w_i\in W_i$. For these insertions, the shifted pseudo-traces of products of
geometrically modified logarithmic $B$-intertwining operators
extend analytically to genus-one $n$-point functions. As the
admissible pseudo-trace data vary, these functions span a space
closed under modular transformations \cite[Theorem~5.5]{HuangLog}.
If $V$ is strongly rational, the same conclusion holds with
ordinary traces.
\end{corollary}
\begin{proof}
By Corollary~\ref{cor:c2}, $B$ is $C_2$-cofinite and has no
negative conformal weights. The first assertion is Huang's
modular invariance theorem \cite[Theorem~5.5]{HuangLog} applied
to $B$. If $V$ is strongly rational, so is $B$ by
Theorem~\ref{thm:rational}, and the semisimple assertion
follows from \cite[Theorem~7.3]{HuangGenusOne}.
\end{proof}

\appendix
\section{Formal products and vacuum coefficients}\label{app:voa}

We verify the expansion conventions used for partial translations
and collisions in Sections~\ref{sec:fields} and~\ref{sec:cutoff}.
All calculations of field products here take place in $V$.
The last argument gives the vacuum partners used to remove the
auxiliary pairs.

\subsection{Regularity}\label{app:regularity}

\begin{proof}[Proof of Lemma~\ref{lem:regularity}]
For homogeneous $a_i,a_j$, the grading gives
$(a_i)_\ell a_j=0$ for $\ell\ge h_i+h_j$. The commutator formula
therefore makes $(z_i-z_j)^{h_i+h_j}$ a locality factor.
In the regularized product, move $Y(a_i,z_i)$ next to the vacuum.
Creation gives $e^{z_iT}a_i$, so no negative powers of $z_i$ occur.
This holds for every $i$, proving joint regularity.
Before multiplication by $\Delta$, each coefficient of total degree
$s$ has weight $\sum_i h_i+s$. Homogeneity of $\Delta$ gives
\eqref{eq:degreeweight}. At fixed weight only one total degree is
possible, so the projection of the power series is a polynomial.
\end{proof}

\subsection{Reexpansion at a common coordinate}\label{app:iterate}

\begin{proof}[Proof of Lemma~\ref{lem:iterate}]
We first establish the identity for an arbitrary fixed homogeneous
vector $v\in V$. We then apply it coefficientwise to the expansion
$W_{I^c}(\mathbf{z}_{I^c})
 =\sum_{\boldsymbol{\gamma}\ge\mathbf{0}}
 v_{\boldsymbol{\gamma}}\mathbf{z}_{I^c}^{\boldsymbol{\gamma}}$.

\smallskip\noindent
\emph{A fixed complementary vector.} For this calculation only, relabel the indices
of $I$ as $1,\ldots,m$, preserving their order, and put $x_m=0$.
Write $\Delta_I$ for the corresponding internal polynomial and
$C(\mathbf{x})$ for $C_I(\mathbf{x})$ in these labels.
Fix $v\in V_d$ and put
\[
 F_v(\mathbf{z})=\Delta_I(\mathbf{z})
       Y(a_1,z_1)\cdots Y(a_m,z_m)v,
 \qquad
 K_i=h_i+d,\qquad
 \Pi(\mathbf{z})=\prod_{i=1}^m z_i^{K_i}.
\]
Insert $Y(v,s)$ immediately before the vacuum and put
\[
 \Psi(\mathbf{z};s)=
 \Delta_I(\mathbf{z})\prod_{i=1}^m(z_i-s)^{K_i}
       Y(a_1,z_1)\cdots Y(a_m,z_m)Y(v,s)\ket{0}.
\]
By Lemma~\ref{lem:regularity}, $\Psi\in V[[\mathbf{z},s]]$.
Setting $s=0$ and using creation gives
$\Psi(\mathbf{z};0)=\Pi(\mathbf{z})F_v(\mathbf{z})$, so
\[
 F_v\in\mathcal A:=
 V[[z_1,\ldots,z_m]][z_1^{-1},\ldots,z_m^{-1}].
\]

There is a linear reexpansion map
\[
 \sigma:\mathcal A\longrightarrow V((t))[[x_1,\ldots,x_{m-1}]]
\]
defined by ordinary substitution $z_i=t+x_i$, $z_m=t$ on
power-series numerators, and by
\[
 z_i^{-1}\longmapsto
 (t+x_i)^{-1}=\sum_{j\ge0}(-1)^j t^{-j-1}x_i^j
 \quad(i<m),\qquad z_m^{-1}\longmapsto t^{-1}.
\]
The numerator substitution preserves total degree, so each joint
coefficient is a finite sum. If a denominator is
$\prod_i z_i^{b_i}$, then the coefficient of
$\mathbf{x}^{\boldsymbol{\alpha}}$ after substitution has no
$t$-power below $-\sum_i b_i-|\boldsymbol{\alpha}|$.
The transformed coordinates are units in $\C((t))[[\mathbf{x}]]$,
so $\sigma$ is independent of the denominator and commutes with
multiplication by scalar polynomials. In particular
$\Pi(t+\mathbf{x},t)\sigma(F_v)=\Psi(t+\mathbf{x},t;0)$.

All scalar factors of $\Psi$ are coordinate differences.
Translation covariance and $T\ket{0}=0$ therefore give
\[
 \left(\sum_{i=1}^m\partial_{z_i}+\partial_s\right)\Psi=T\Psi,
 \qquad
 \Psi(\mathbf{z}+u;s+u)=e^{uT}\Psi(\mathbf{z};s).
\]
The second equality is the formal Taylor identity for the first.
Both sides are joint power series. Substituting
$z_i=x_i$ for $i<m$, $z_m=0$, $s=-t$, and $u=t$ consequently gives
\[
 \Psi(t+\mathbf{x},t;0)=e^{tT}\Psi(\mathbf{x},0;-t).
\]
In $\Psi$, locality allows $Y(v,s)$ to be moved to the far left
while retaining every factor $(z_i-s)^{K_i}$ in its original
orientation. The field $Y(a_m,z_m)$ is then rightmost, so
constant-coefficient extraction in $z_m$ uses creation and yields
\[
 \Psi(\mathbf{x},0;-t)=
 \Pi(t+\mathbf{x},t)Y(v,-t)C(\mathbf{x}).
\]
For each coefficient $c_{\boldsymbol{\alpha}}$ of $C$, skew symmetry
gives $e^{tT}Y(v,-t)c_{\boldsymbol{\alpha}}
=Y(c_{\boldsymbol{\alpha}},t)v$. Lower truncation makes the exponential
action well defined at each coefficient. Combining the identities gives
\[
 \Psi(t+\mathbf{x},t;0)=
 \Pi(t+\mathbf{x},t)Y(C(\mathbf{x}),t)v.
\]
The constant coefficient of $\Pi(t+\mathbf{x},t)$ in $\mathbf{x}$ is
$t^{\sum_iK_i}$, so $\Pi(t+\mathbf{x},t)$ is a unit in
$\C((t))[[\mathbf{x}]]$. Cancellation gives
\[
 \sigma(F_v)=Y(C(\mathbf{x}),t)v.
\]
For the two-field calculation with these expansion conventions, see
\cite[correction to Theorem~3.2.1]{FBZcorr}.

\smallskip\noindent
\emph{Reassembling the complementary product.}
Restore the original labels and write $\mathbf{y}=\mathbf{z}_{I^c}$.
Locality, with the original orientations of all factors retained,
gives
\[
 H_V(\mathbf{z})(e)=
 \Delta_{\rm cross}(\mathbf{z}_I,\mathbf{y})
 \sum_{\boldsymbol{\gamma}\ge\mathbf{0}}
 F_{v_{\boldsymbol{\gamma}}}(\mathbf{z}_I)
       \mathbf{y}^{\boldsymbol{\gamma}}.
\]
Each $v_{\boldsymbol{\gamma}}$ is homogeneous when nonzero.
Write the scalar cross factor as the finite sum
\[
 \Delta_{\rm cross}(\mathbf{z}_I,\mathbf{y})
 =\sum_{\boldsymbol{\delta}\in D}
      p_{\boldsymbol{\delta}}(\mathbf{z}_I)
      \mathbf{y}^{\boldsymbol{\delta}}.
\]
For any fixed $\boldsymbol{\beta}$, coefficient extraction gives
\[
 [\mathbf{y}^{\boldsymbol{\beta}}]H_V(\mathbf{z})(e)
 =
 \sum_{\substack{\boldsymbol{\delta}\in D\\
                 \boldsymbol{\delta}\le\boldsymbol{\beta}}}
 p_{\boldsymbol{\delta}}(\mathbf{z}_I)
 F_{v_{\boldsymbol{\beta}-\boldsymbol{\delta}}}(\mathbf{z}_I).
\]
Apply $\sigma$ to this finite sum. On the left it is substitution
in the power series $H_V$. On each term on the right it commutes
with the scalar polynomial and gives the iterate already computed.
This proves \eqref{eq:iterateformula} at each
$\mathbf{y}^{\boldsymbol{\beta}}$, allowing the coordinate-pole bounds
to depend on the complementary coefficient.
\end{proof}

\subsection{Vacuum coefficients}\label{app:vacuum}

\begin{proof}[Proof of Lemma~\ref{lem:detect}]
Simplicity enters through the ideal of states invisible to all vacuum
coefficients. Let $\pi_0:V\to V_0=\C\ket{0}$ be the graded projection, and put
\[
 J=\{v\in V:\pi_0Y(a,z)v=0\text{ for all }a\in V\}.
\]
For homogeneous $a$, the contributions of distinct homogeneous
components of $v$ to $\pi_0Y(a,z)v$ occur at distinct powers of
$z$. Thus $J$ is graded. For $a,b\in V$ and $v\in J$,
the definition of $J$ gives
\[
 \pi_0Y(Y(a,x)b,t)v=0.
\]
By formal associativity \cite[Proposition 2.6]{Frenkel}, this
iterate and $\pi_0Y(a,z_1)Y(b,z_2)v$ are expansions of the same
element of
\[
 V_0[[z_1,z_2]]
 [z_1^{-1},z_2^{-1},(z_1-z_2)^{-1}].
\]
For the iterate, set $z_1=t+x$, $z_2=t$ and expand
$(t+x)^{-1}$ in nonnegative powers of $x/t$.
Clearing the displayed denominators reduces injectivity of the
expansion to the invertible linear change of variables on the
power-series numerator. Thus
$\pi_0Y(a,z_1)Y(b,z_2)v=0$.
Taking coefficients shows that $b_mv\in J$ for all $b\in V$ and $m\in\mathbb Z$.
Since $\pi_0T=0$, translation covariance gives
\[
 \pi_0Y(a,z)Tv=-\partial_z\pi_0Y(a,z)v=0.
\]
Thus $J$ is stable under all left modes and $T$, hence is an ideal
by skew symmetry. Since $\ket{0}\notin J$, simplicity gives $J=0$.

Since $u\notin J$, some homogeneous $a\in V_k$ satisfies
$a_{h+k-1}u\in\C^\times\ket{0}$.
Necessarily $k>0$, since $Y(\ket{0},z)u=u$. Skew symmetry
$Y(u,z)a=e^{zT}Y(a,-z)u$ gives
\[
 u_{h+k-1}a=(-1)^{h+k}a_{h+k-1}u.
\]
The terms involving positive powers of $T$ have zero vacuum component.
Take $w=a$.
\end{proof}

\begingroup
\small
\raggedright

\endgroup
\end{document}